\documentclass[10pt]{article}
\usepackage{definitions}

\newcommand{\lc}{\ell_c}
\newcommand{\uc}{u_c}

\newcommand{\R}{{\mathbb R}}

\SetKwInOut{Input}{Input}
\SetKwComment{Comment}{$\triangleright$\ }{}

\title{MPAX: Mathematical Programming in JAX}
\author{Haihao Lu\thanks{MIT, Sloan School of Management (haihao@mit.edu).} \and Zedong Peng\thanks{MIT, Sloan School of Management (zdpeng@mit.edu).} \and Jinwen Yang\thanks{University of Chicago, Department of Statistics (jinweny@uchicago.edu).}}
\date{}

\begin{document}
\maketitle

\begin{abstract}
We present MPAX (Mathematical Programming in JAX), an open-source first-order solver for large-scale linear programming (LP) and convex quadratic programming (QP) built natively in JAX. The primary goal of MPAX is to exploit modern machine learning infrastructure for large-scale mathematical programming, while also providing advanced mathematical programming algorithms that are easy to integrate into machine learning workflows.
MPAX implements two PDHG variants, $\mathrm{r^2HPDHG}$ for LP and rAPDHG for QP, together with diagonal preconditioning, adaptive restarts, adaptive step sizes, primal-weight updates, infeasibility detection, and feasibility polishing. Leveraging JAX’s compilation and parallelization ecosystem, MPAX provides across-hardware portability, batched solving, distributed optimization, and automatic differentiation.
We evaluate MPAX on CPUs, NVIDIA GPUs, and Google TPUs, observing substantial GPU speedups over CPU baselines and competitive performance relative to GPU-based codebases on standard LP/QP benchmarks. Our numerical experiments further demonstrate MPAX's capabilities in high-throughput batched solving, near-linear multi-GPU scaling for dense LPs, and efficient end-to-end differentiable training. The solver is publicly available at \url{https://github.com/MIT-Lu-Lab/MPAX}.
\end{abstract}

\section{Introduction}

Modern machine learning and deep learning have advanced at an extraordinary pace, driven not only by algorithmic innovations but also by the maturation of accelerator-centric software and hardware infrastructure~\cite{paszke2019pytorch}. Today’s ML systems operate within tightly integrated pipelines built on XLA-based compilation, batching and vectorization primitives, multi-device distribution, and automatic differentiation. This ecosystem has enabled unprecedented scalability, throughput, and ease of deployment~\cite{openxla_google_blog_2023}.

In contrast, traditional mathematical optimization has historically been engineered as a standalone computational stage, often residing outside the ML execution stack~\cite{ning2019optimization,freitas2009prediction}. Classical solvers rely heavily on CPU-oriented designs, sparse linear algebra factorizations, and host-driven control flow. As a result, integrating optimization into modern accelerator pipelines typically requires shuttling data between device-resident ML models and host-based solvers, interrupting the otherwise unified computation flow. In many applied systems, this leads to a pattern where an ML model predicts inputs (e.g., costs, revenues, or consumptions), transfers them from the GPU back to the CPU, invokes an external solver such as Gurobi~\cite{gurobi}, and then passes the solution back into the ML pipeline. Although this approach is workable, it has certain drawbacks: it incurs nontrivial data-movement overhead, introduces latency, breaks differentiability, and prevents native batching or parallel execution. 

To address this, we introduce MPAX (Mathematical Programming in JAX), a JAX-native {first-order} solver for linear programming (LP) and quadratic programming (QP). MPAX is designed with a core objective:
\begin{center}
    \textit{To develop an efficient optimization framework that leverages modern ML infrastructure and its advanced capabilities to facilitate diverse application scenarios.}
\end{center}

At its algorithmic core, MPAX is designed as a purely first-order solver, relying exclusively on gradient information rather than Hessian information or matrix factorizations. The LP component adopts the algorithmic framework of cuPDLPx~\cite{lu2025cupdlpx}, and the QP component builds on PDQP~\cite{lu2025practical}. This design choice is motivated by several factors. First, recent advances in first-order methods for LPs and QPs have demonstrated that performance on large-scale problems can be on par with traditional solvers based on simplex and interior-point methods. Second, first-order kernels map naturally to modern accelerator architectures; operations such as matrix-vector products and element-wise updates are highly parallelizable, allowing MPAX to saturate the compute throughput of GPUs and TPUs. Finally, this approach provides a lightweight memory footprint, enabling the solution of massive-scale optimization instances that would otherwise exceed the memory constraints of dedicated hardware.

MPAX is explicitly designed to leverage core transformations in the JAX ecosystem, including just-in-time compilation, vectorization, distributed optimization, and automatic differentiation. When applied to LP/QP solving, these transformations translate into several advanced capabilities that are difficult (at least without substantial engineering effort) to realize in traditional solver implementations.
\begin{itemize}
    \item \textbf{Across-hardware portability and hardware acceleration.} MPAX runs natively on CPUs, GPUs, and TPUs, utilizing XLA-based just-in-time (JIT) compilation to lower high-level Python updates into fused, hardware-optimized executables. Our empirical evaluation demonstrates that MPAX achieves a significant improvement in solve time on modern GPUs compared with CPU execution.
    \item \textbf{Batched solving.} MPAX supports batched solving. It can process large batches of same-shape LP/QP instances in parallel through automatic vectorization, allowing the solver to saturate accelerator throughput and avoid host-side loops. In a two-stage stochastic LP experiment, MPAX exhibits much more favorable scaling with batch size than a CPU-based solver, yielding substantial speedups as the batch dimension grows.
    \item \textbf{Distributed optimization.} MPAX supports distributed execution through Single-Program Multiple-Data (SPMD) computation with automatic data sharding across a device mesh. This enables both memory scaling, allowing larger LP/QP instances to fit across devices, and performance scaling. In a large-scale dense LP experiment, we observe near-linear speedup as the number of hardware devices increases.
    \item \textbf{Auto-differentiation.} By unrolling iterations or using surrogate differentiation schemes, MPAX can compute or approximate derivatives of the optimal solution with respect to input parameters. This makes it possible to embed optimization directly within differentiable pipelines, enabling modern ML-driven decision-making frameworks.
\end{itemize}

Lastly, we also benchmark MPAX against cuPDLPx~\cite{lu2025cupdlpx} and PDQP~\cite{lu2025practical} in standalone (single-instance) solves. Overall, MPAX is 1.05×–1.18× slower than cuPDLPx in runtime, while outperforming PDQP by 1.25×–1.6×. These results indicate that MPAX remains competitive even in traditional solver-centric settings, while additionally offering additional JAX-native capabilities beyond traditional solvers.


\subsection{Related works}

Traditional linear and quadratic programming solvers, essentially based on the simplex method and interior-point methods (IPMs), are predominantly CPU-based. These methods are widely used in industrial solvers, such as Gurobi~\cite{gurobi}, Mosek~\cite{mosek}, COPT~\cite{ge2022cardinal} and FICO Xpress~\cite{berthold2018parallelization}, due to their robustness and mature implementations. However, extending these methods to modern accelerators is often non-trivial when performance is dominated by sparse linear system solves and matrix factorizations, which are less amenable to fine-grained parallelism and frequently involve irregular memory access patterns.

Motivated by large-scale instances where factorization-based methods can become the computational bottleneck, first-order methods have received renewed attention for LP/QP. A representative class is the ADMM-based approaches, where a matrix factorization is performed once and reused across iterations. Examples include OSQP~\cite{osqp,osqp-gpu}, SCS~\cite{o2016conic}, ABIP~\cite{deng2025enhanced} and other ADMM-style implementations. Another prominent class is factorization-free primal--dual first-order methods, which are typically dominated by sparse matrix--vector multiplications and element-wise operations. This direction has led to several GPU-native solvers for large-scale LP and QP, including cuPDLP.jl~\cite{lu2023cupdlp,lu2023cupdlpc}, cuPDLPx~\cite{lu2025cupdlpx}, and HPR-LP~\cite{chen2024hpr} for LP, as well as PDQP.jl~\cite{lu2025practical}, PDHCG~\cite{huang2025restarted}, and HPR-QP~\cite{chen2025hpr} for QP. In addition, as a new algorithmic choice, variants of PDLP have been integrated into several industrial solver stacks, including NVIDIA cuOpt, COPT, Gurobi, FICO Xpress, and Knitro, perhaps among others. Most of these systems are implemented in performance-oriented languages such as C/C++ or Julia and are released as standalone solvers or language-specific libraries. 

In parallel to solver implementations that expose GPU acceleration through dedicated codebases, another thread builds optimization primitives directly within machine learning programming frameworks, such as JAX and PyTorch~\cite{paszke2019pytorch}. For example, qpth~\cite{amos2017optnet}, osqpth, and qpax~\cite{tracy2024differentiability} provide QP-related primitives, while general toolkits such as JAXopt~\cite{blondel2022efficient} and Optax~\cite{deepmind2020jax} offer framework-native optimization building blocks. These toolkits primarily target optimization for machine learning and deep learning objectives, rather than specialized large-scale mathematical programming solvers, and their implementations that are not targeted for challenging LP/QP instances. Compared with standalone solver libraries, framework-native implementations can directly operate on device-resident arrays and compose with compilation, vectorization/batching, and multi-device execution mechanisms provided by the underlying frameworks. MPAX belongs to this line of work and implements advanced first-order primal-dual algorithms for general-purpose LP/QP.

\section{Preliminary: Primal-Dual First-Order Algorithms for LP and QP}\label{sec:preliminary}

MPAX currently supports solving general linear programs and quadratic programs in the following forms
\begin{equation}\label{eq:qp}
    \begin{aligned}[c]
    \min_{x}~~ &~ \frac{1}{2} x^\top Q x + c^\top x \\
\text{s.t.}~~ &~ l_c \leq Ax \leq u_c \\
&~ l_v \leq x \leq u_v  ,
    \end{aligned}
\end{equation}
where $Q \in \mathbb R^{n \times n}$, $A \in \mathbb R^{m \times n}$, $c \in \mathbb R^{n}$, $l_c \in (\mathbb R \cup \{ -\infty \})^{n}$, $u_c \in (\mathbb R \cup \{ \infty \})^{n}$, $l_v \in (\mathbb R \cup \{ -\infty \})^{n}$, $u_v \in (\mathbb R \cup \{ \infty \})^{n}$. When $Q=0$, \eqref{eq:qp} reduces to a linear program.
MPAX solves these LP and QP by addressing their primal-dual form:

\begin{flalign}\label{eq:primal-dual-qp}
\min_{x \in X} \max_{y \in Y}\ L(x,y) := \frac{1}{2} x^\top Q x + c^\top x - y^\top A x - p(y;l_c,u_c)  \ ,
\end{flalign}

where $p(y;-u_c,-l_c) = u_c^\top y^- -l_c^\top y^+ $, $X := \{x \in \mathbb R^n : l_v \leq x \leq u_v \}$, and 
\begin{equation*}
\begin{aligned}
    &Y_i := \begin{cases}
\{ 0 \} & (\lc)_i = -\infty, ~ (\uc)_i = \infty,  \\
\R^{-} & (\lc)_i = -\infty, ~ (\uc)_i \in \R ,\\
\R^{+} & (\lc)_i \in \R, ~ (\uc)_i = \infty, \\
\R & \text{otherwise}.
\end{cases}
\end{aligned}
\end{equation*}

\subsection{Restarted Reflected Halpern PDHG for LP}\label{sec:lp-pdhg}

MPAX implements the restarted Halpern PDHG with reflection ($\mathrm{r^2HPDHG}$)~\cite{lu2024restarted} to solve LP problems, based on the implementation developed in cuPDLPx~\cite{lu2025cupdlpx}. The full algorithmic details and engineering considerations are provided in \cite{lu2025cupdlpx}, while the intuition and theoretical guarantees of $\mathrm{r^2HPDHG}$ are presented in \cite{lu2024restarted}. In this subsection, we present a high-level overview of the algorithmic structure (see \cite{lu2025cupdlpx,lu2024restarted} for comprehensive details).

For LP, the quadratic term $Q=0$ in \eqref{eq:qp}, and the Primal Dual Hybrid Gradient (PDHG) update rule  is as follows:
\begin{equation}\label{eq:pdhg-lp} 
    \begin{aligned}
        & x^{k+1}\leftarrow \text{proj}_{X}(x^k-\tau (c-A^\top y^k)) \\ 
        & y^{k+1}\leftarrow y^k-\sigma A(2x^{k+1}-x^k)-\sigma\operatorname{proj}_{[-u_c,-l_c]}\pran{\sigma ^{-1}y^k-A(2x^{k+1}-x^k)}\ , 
    \end{aligned}
\end{equation}
where $\tau$ and $\sigma$ are the primal step size and dual step size, respectively.
However, the performance of the vanilla PDHG method for solving LPs is insufficient for general-purpose solvers~\cite{applegate2021practical}. To address this and motivated by theoretical understandings~\cite{applegate2023faster,applegate2024infeasibility}, several algorithmic enhancements have been proposed, including a restart scheme, adaptive step sizes, primal weighting, and preconditioning, which significantly improve the numerical performance of the algorithm, together with GPU acceleration~\cite{applegate2021practical,lu2023cupdlp,lu2023cupdlpc}. 

More recently, a variant of the algorithm, Halpern PDHG with reflection ($\mathrm{r^2HPDHG}$), has been introduced~\cite{lu2024restarted,chen2024hpr}. The update rule for Halpern PDHG with reflection is given by:
\begin{equation}\label{eq:hrpdhg}
    (x^{k+1}, y^{k+1})\leftarrow  \frac{k+1}{k+2}( (1 + \gamma) \text{PDHG}((x^k, y^k))- \gamma (x^k, y^k))+\frac{1}{k+2}(x^0, y^0) \ ,
\end{equation}
where $\gamma \in [0, 1]$ is the reflection coefficient (it is also called the over-relaxation parameter in variational inequality literature) and $PDHG(\cdot)$ represent the update rule defined in \eqref{eq:pdhg-lp}.
In this formulation, the next iterate is a weighted average of a reflected PDHG step ($2 \cdot \text{PDHG}(z^k)-z^k$) and the initial solution $z^0$. 
The $\mathrm{r^2}$HPDHG algorithm is summarized in Algorithm \ref{alg:r2hpdhg}. The method adopts a nested-loop structure: the inner loop iteratively performs the PDHG update followed by the Halpern update with reflection, while the outer loop governs the restart scheme. When the restart criterion is met, the inner loop terminates, and the initial point is reset to the current iterate to start the next cycle. The algorithm stops once the termination criterion is satisfied, typically based on the KKT residual. 

\setlength{\algomargin}{1em}
\begin{algorithm}[H]
\caption{Reflected restarted Halpern PDHG ($\mathrm{r^2HPDHG}$)}
\label{alg:r2hpdhg}
\DontPrintSemicolon

\Input{Linear program \eqref{eq:qp}, initial point $(x^{0,0},y^{0,0})$,  termination tolerance $\epsilon$.}
\For{$n \gets 0, 1, \dots$}{
    $k \gets 0$ \Comment*[r]{Inner loop counter}
    \Repeat{restart condition holds}{
        $\hat x^{n,k+1} \gets \operatorname{proj}_{X}(x^{n,k}-\tau (c-A^\top y^{n,k}))$ \Comment*[r]{PDHG update}
        $\hat y^{n,k+1} \gets y^{n,k}-\sigma A(2x^{n,k+1}-x^{n,k}) - \sigma \operatorname{proj}_{[-u_c,-l_c]}\left(\sigma ^{-1}y^{n,k}-A(2x^{n,k+1}-x^{n,k})\right)$\;
        \If{$\text{KKT}(\hat x^{n,k+1}, \hat y^{n,k+1}) \le \epsilon$} 
        {
            \Return{$(\hat x^{n,k+1},\hat y^{n,k+1})$} \Comment*[r]{Check termination}
        }
        $x^{n,k+1} \gets \frac{k+1}{k+2}(2 \hat x^{n,k+1}-x^{n,k})+\frac{1}{k+2}x^{n,0}$ \Comment*[r]{Halpern update \& reflection}
        $y^{n,k+1} \gets \frac{k+1}{k+2}(2 \hat y^{n,k+1}-y^{n,k})+\frac{1}{k+2}y^{n,0}$\;
        
        $k \gets k+1$\;
    }
    $(x^{n+1,0}, y^{n+1,0}) \gets (\hat x^{n,k}, \hat y^{n,k})$ \Comment*[r]{Restart}

}
\end{algorithm}

\subsection{Restarted Accelerated PDHG for QP}\label{sec:qp-pdhg}

For convex quadratic programming, MPAX implements the restarted accelerated PDHG algorithm (rAPDHG) described in \cite{lu2025practical}. Below we provide a high-level overview of rAPDHG; for a full algorithmic description, implementation details, and theoretical motivation, we refer the reader to \cite{lu2025practical}.

For a QP in the form of \eqref{eq:qp}, a single PDHG iteration is defined as follows:

\begin{equation}\label{eq:pdhg-qp} 
    \begin{aligned}
        & x^{k+1}\leftarrow \text{proj}_{X}(x^k-\tau (Qx^{k} + c-A^\top y^k)) \\ 
        & y^{k+1}\leftarrow y^k-\sigma A(2x^{k+1}-x^k)-\sigma\operatorname{proj}_{[-u_c,-l_c]}\pran{\sigma ^{-1}y^k-A(2x^{k+1}-x^k)}\ , 
    \end{aligned}
\end{equation}
where $\tau$ and $\sigma$ are the primal step size and dual step size, respectively.
As in the LP case, vanilla PDHG is typically not sufficient for a robust and general-purpose QP solver~\cite{lu2025practical}. A restarted accelerated PDHG method has been proposed that attains an optimal linear convergence rate among a broad class of primal--dual methods for convex QPs~\cite{lu2025practical}. The acceleration is introduced through a momentum term in the primal update. In particular, with a running average $\bar x^k$, the momentum-augmented primal update can be written as

\begin{equation}\label{eq:pdhg-momentum-qp}
    x^{k+1}\leftarrow \text{proj}_{X}(x^k-\tau (Q ( (1-\beta) \bar x^{k} + \beta x^{k}) + c-A^\top y^k))\ ,
\end{equation}
where $\beta$ is the momentum coefficient. Beyond momentum, rAPDHG typically combines several additional mechanisms, including restart, averaging, adaptive step sizes, primal weighting, and preconditioning, to improve robustness and performance across heterogeneous QP instances~\cite{lu2025practical}. Algorithm~\ref{alg:rapdhg} presents the restarted accelerated PDHG (rAPDHG) algorithm implemented in MPAX. The method follows a nested-loop structure. In the inner loop, the algorithm repeatedly applies PDHG updates together with a weighted-averaging step, while incorporating momentum acceleration on the primal iterates. The outer loop implements the restart mechanism: once the restart condition is triggered, the inner loop stops and the next cycle is initialized from the current weighted-average iterate. The procedure terminates when the stopping criterion is met, typically evaluated via the KKT residual. Further algorithmic details and theoretical guarantees are provided in \cite{lu2025practical}.

\begin{algorithm}[ht!]
        \SetAlgoLined
        {\bf Input:} Quadratic program \eqref{eq:qp}, initial point $(x^{0,0},y^{0,0})$, step-size $\{(\beta_k,\theta_k,\tau_k,\sigma_k)\}$\;
            \For{$n \gets 0, 1, \dots$}{
                $k \gets 0$ \Comment*[r]{Inner loop counter}
            \Repeat{restart condition holds}{
                $x^{n,k+1}\leftarrow \text{proj}_{X}(x^{n,k}-\tau_k (Q ( (1-\beta_k) \bar x^{n,k} + \beta_k x^{n,k}) + c-A^\top y^{n,k}))$ \Comment*[r]{PDHG update \& Momentum}
                \vspace{0.075cm}
                $y^{n,k+1}\leftarrow y^{n,k}-\sigma_k A(2x^{n,k+1}-x^{n,k})-\sigma_k\operatorname{proj}_{[-u_c,-l_c]}\pran{\sigma_k ^{-1}y^{n,k}-A(2x^{n,k+1}-x^{n,k})}$\;
                \vspace{0.075cm}
                \vspace{0.075cm}
                $\bar x^{n,k+1}=(1-\theta_k)\bar x^{n,k}+\theta_k x^{n,k+1}$ \Comment*[r]{Weighted Average}
                \vspace{0.075cm}
                $\bar y^{n,k+1}=(1-\theta_k)\bar y^{n,k}+\theta_k y^{n,k+1}$\;
                \If{$\text{KKT}(\bar x^{n,k+1}, \bar y^{n,k+1}) \le \epsilon$} 
                {
                    \Return{$(\bar x^{n,k+1},\bar y^{n,k+1})$} \Comment*[r]{Check termination}
                }
                $k \gets k+1$\;
                
            }
        $(x^{n+1,0}, y^{n+1,0})\gets (\bar x^{n,k},\bar y^{n,k})$  \Comment*[r]{Restart}
        $(\bar x^{n+1,0},\bar y^{n+1,0}) \gets (\bar x^{n,k},\bar y^{n,k})$\;
        }
        \caption{Restarted accelerated PDHG (rAPDHG)}
        \label{alg:rapdhg}
\end{algorithm}

\section{MPAX: Math Programming in JAX}\label{sec:mpax}

\href{https://github.com/MIT-Lu-Lab/MPAX}{MPAX} is a hardware-accelerated, batchable, distributable, and differentiable solver built entirely in JAX, designed to address classic optimization problems encountered in data science and machine learning.
MPAX has implemented $\mathrm{r^2}$HPDHG algorithm for LP and rAPDHG algorithm for QP as described in Section \ref{sec:preliminary}.

\begin{figure}[bht]
\begin{center}
\centerline{\includegraphics[width=0.65\columnwidth]{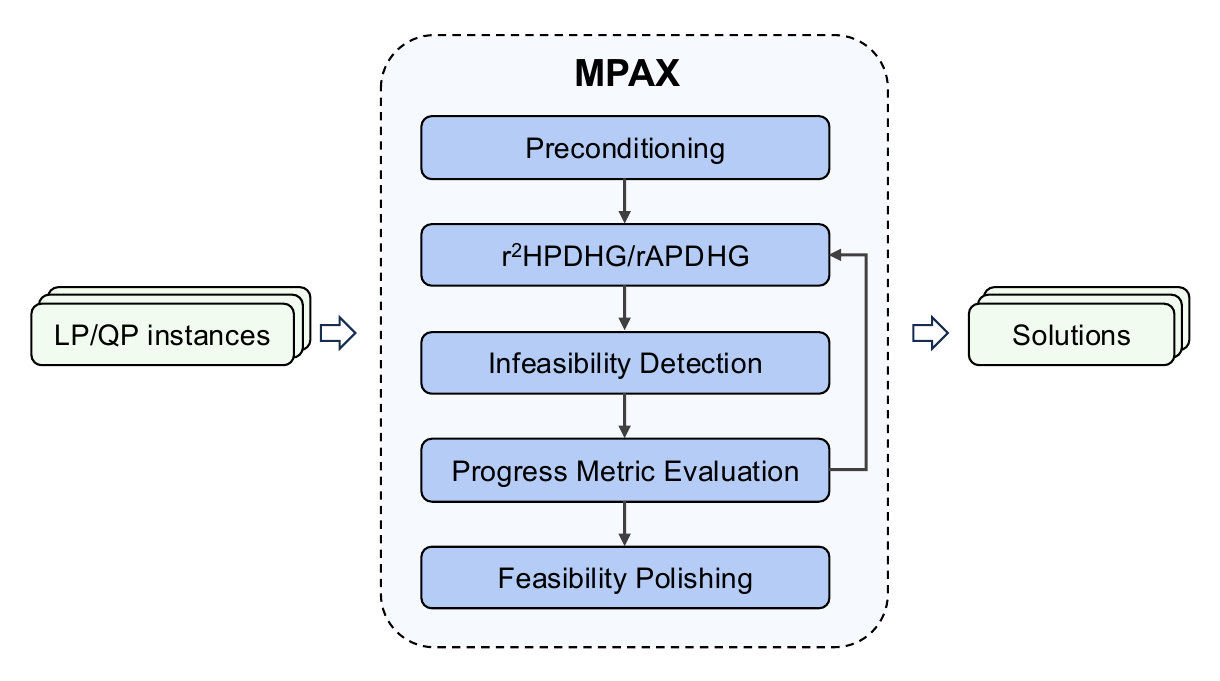}}
\caption{Design of MPAX}
\label{fig:mpax design} 
\end{center}
\end{figure}

The overall design of MPAX is illustrated in Figure \ref{fig:mpax design}. The process begins with preconditioning to improve the condition number of the input LP/QP problem. By default, MPAX applies two diagonal preconditioners: Ruiz scaling \cite{ruiz2001scaling} and Pock and Chambolle’s diagonal scaling methods \cite{pock2011diagonal}.
After preconditioning, MPAX runs $\mathrm{r^2}$HPDHG iterations on the scaled LP or rAPDHG iterations on the scaled QP. During the iterations, three key conditions are periodically evaluated: termination, restart, and infeasibility detection. Iterations terminate when the relative KKT error, comprising primal feasibility, dual feasibility, and duality gap, meets the specified tolerance. These algorithmic mechanisms follow the designs used in cuPDLPx~\cite{lu2025cupdlpx} and PDQP~\cite{lu2025practical}, including step-size rules, restart conditions, and primal-weight updates. 

One notable feature of MPAX that differs from cuPDLPx~\cite{lu2025cupdlpx} and PDQP~\cite{lu2025practical} is the inclusion of a feasibility polishing step. In machine learning workflows, solutions of moderate optimality accuracy are often acceptable, but feasibility remains crucial. That is, feasibility tolerances typically need to be tighter than optimality tolerances. To address this, feasibility polishing, first introduced in~\cite{applegate2025pdlp} in the context of PDLP, is employed as a post-processing step to improve feasibility at the cost of potentially relaxing optimality. This approach exploits two observations: (1) primal–dual first-order methods converge more rapidly on feasibility problems than on optimality problems, and (2) starting PDLP from a reasonable initial point leads to convergence to a nearby optimal solution. Feasibility polishing therefore solves the feasibility-only problem (objective set to zero) using PDLP, allowing the iterations to focus entirely on reducing constraint violation. MPAX adopts a polishing procedure aligned with the description in~\cite{applegate2025pdlp}.

The remainder of this section introduces the advanced features of MPAX, including across-hardware portability and hardware acceleration, batched solving, distributed optimization, automatic differentiation, and additional utilities. Most of these capabilities are not supported by previous Julia- or C-based implementations~\cite{lu2023cupdlp,lu2024restarted,lu2025cupdlpx}. A detailed description of solver options is presented in Appendix \ref{appendix:solver options}.

\subsection{Across-hardware portability and hardware acceleration}
Built on JAX, MPAX natively runs on multiple types of hardware, including CPUs, GPUs and TPUs. It takes full advantage of the just-in-time (JIT) compilation powered by the XLA compiler, which lowers high-level Python programs to highly optimized, hardware-specific executables for CPUs, GPUs, and TPUs. This compilation pipeline enables operation fusion, parallel execution, and target-specific optimizations, yielding substantial speedups. Crucially, since first-order methods are often memory-bound, XLA's operation fusion significantly reduces memory bandwidth usage by keeping intermediate results in on-chip registers rather than round-tripping to main memory. On NVIDIA GPUs, for instance, this results in sequences of element-wise operations being compiled into efficient, fused CUDA kernels.

This compilation paradigm is particularly effective for $\mathrm{r^2}$HPDHG and rAPDHG. As first-order methods, their update rules depend only on explicit gradient computations and simple proximal/projection operations, making them well-suited for efficient implementation in JAX. Beyond PDHG updates, MPAX further implements algorithmic heuristics, including restart strategies, adaptive step sizes, and line search, utilizing JAX control-flow primitives. This allows the entire solver loop to be compiled into a single computational graph, thereby eliminating the overhead of kernel launching and host-device synchronization associated with Python-based iteration. Overall, the implementations of $\mathrm{r^2}$HPDHG and rAPDHG in MPAX are engineered to maximize the benefits of JAX's JIT compilation for hardware acceleration.

\subsection{Batched solving}
MPAX supports efficient batch-solving of LP/QP instances with identical shapes on a single accelerator, such as a GPU or TPU, through JAX's automatic vectorization transformation (\texttt{vmap}). When multiple LP/QP instances must be solved simultaneously, a traditional solver would process them sequentially, handling one instance at a time. In contrast, MPAX performs true batch solving by stacking problem data along a leading batch dimension and relying on JAX’s batching rules to lift each primitive operation into its batched equivalent. For example, multiple vector-vector multiplications can be transformed into a matrix-vector multiplication, allowing multiple operations to be executed in parallel within a single batched computation. 
This capability is particularly critical for GPUs and TPUs, as solving a single small-to-medium scale instance often fails to saturate the device's massive computational capacity. By processing batches, MPAX maximizes hardware occupancy and throughput, and eliminates the need for explicit for-loops, reducing Python interpreter overhead. Consequently, it unlocks the capability for massive parallelism in advanced applications, such as decomposition methods for block-angular structure instances~\cite{rahmaniani2017benders}, batch training in end-to-end decision making~\cite{amos2017optnet}, strong branching in mixed-integer programming~\cite{dey2024theoretical}, and optimization-based bound tightening in global optimization~\cite{gleixner2017three}.

\subsection{Distributed optimization}
MPAX also supports distributed optimization (i.e., execution across multiple GPUs) by combining just-in-time compilation with data-sharding techniques under the Single Program, Multiple Data (SPMD) framework. Traditionally, enabling multi-GPU execution requires substantial engineering effort: matrices and vectors must be manually partitioned, distributed across devices, and updated synchronously in a distributed loop. In contrast, MPAX integrates with JAX’s sharding API, enabling multi-device execution with minimal user intervention. This design can accelerate the matrix–vector multiplication operation, which is often bottlenecked by the memory bandwidth of a single accelerator.

Beyond performance benefits, distributed optimization can also address the memory capacity constraint. GPUs and TPUs typically provide far less memory than high-end CPUs. Therefore, distributing data across multiple devices is often essential for fitting large-scale LP instances. By sharding the most memory-intensive components, primarily the constraint matrix and objective matrix, MPAX can handle problem sizes that exceed the memory of a single device.

The multi-device workflow begins by defining a logical mesh of accelerators with named parallel axes. The matrix is then partitioned across these axes, distributing storage and computation evenly. Once sharded, the data is passed into the JIT-compiled solve function. The JAX compiler automatically inserts the required collective communication primitives (e.g., all-reduce) to synchronize partial products and maintain consistency across devices. This approach ensures efficient memory usage, coordinated updates, and high throughput, all while preserving a clean and simple high-level solver interface.

A current limitation is that distributed optimization for sparse matrices requires more specialized sharding strategies for efficiency. Improving sparse-data sharding remains an item for future development.

\subsection{Auto-differentiation}
\textsc{MPAX} supports auto-differentiation to compute or approximate the derivatives of the optimal solution with respect to input parameters~\cite{baydin2018automatic}. This capability is essential for modern end-to-end decision-making pipelines built on \emph{differentiable optimization}~\cite{agrawal2019differentiable,amos2017optnet}. Differentiable optimization embeds mathematical programs directly into learning architectures and treats the solver as part of the computational graph. By propagating gradients through the optimization layer, models can be trained using downstream decision quality rather than relying only on surrogate prediction losses. This provides a principled framework for decision-focused learning and unifies prediction and optimization in a single pipeline.

For QP, \textsc{MPAX} provides explicit differentiation by unrolling the iterative algorithm (e.g., via \texttt{jax.lax.scan}). 
In contrast, LP poses additional challenges. LP solutions can be discontinuous with respect to problem parameters. Small perturbations may cause optimal bases to change abruptly, making the solution mapping non-differentiable. To address this, surrogate decision losses such as Smart Predict-then-Optimize+ (SPO+) loss~\cite{elmachtoub2022smart} and perturbed Fenchel-Young loss~\cite{berthet2020learning} are used to approximate gradients. These surrogate losses are incorporated into the automatic differentiation system using JAX’s custom derivative rules, enabling a differentiable forward–backward workflow: a forward pass computes optimal solutions or decision losses, and a backward pass computes approximate gradients.

\subsection{Miscellaneous}
\begin{itemize}

    \item {\bf Matrix format.} MPAX supports both dense and sparse formats for the objective and constraint matrix, providing flexibility to handle a wide range of LP/QP structures efficiently. The computational bottleneck of both $\mathrm{r^2}$HDPHG and rAPDHG algorithms lies in matrix-vector multiplication, which has been fully optimized in JAX for both sparse and dense matrices to maximize performance. However, mismatches between the matrix format and the actual data structure, such as storing a dense matrix in a sparse format or vice versa, can significantly impact computational efficiency. To address this, MPAX automatically follows the input matrix format to ensure that computations align with the structure of the given LP/QP problem, thereby maintaining optimal efficiency and performance.

    \item {\bf Warm start.} By default, MPAX initializes with all-zero vectors as the starting point. However, it also supports a warm start with user-provided solutions, which may include a primal solution, a dual solution, or a pair of both. In many machine learning applications, solvers are frequently used to solve a series of similar problem instances where the input data undergoes only minor changes. Although the choice of warm-start strategy is typically problem-dependent, providing a well-suited warm start for both rAPDHG and $\mathrm{r^2}$HPDHG can significantly accelerate their convergence.
    \item {\bf Determinism.} In JAX, floating-point computations on accelerators may produce non-deterministic results due to factors such as XLA auto-tuning, GPU/TPU atomics, and scatter operations. Additionally, \texttt{jax.jit} can introduce variability as a result of optimizations performed by the XLA compiler. Setting \texttt{xla\_gpu\_deterministic\_ops=true} can ensure deterministic results in JAX.
    While these differences are typically minor for each iteration, they can accumulate over time, leading to noticeable variations in convergence.

    \item {\bf Precision.} MPAX uses single-precision (FP32) by default, consistent with JAX's standard setting. For more robust numerical performance and convergence, double-precision (FP64) can be enabled by setting \texttt{jax\_enable\_x64} to True.

    \item \textbf{External access.} With growing interest from the machine learning community, MPAX is now integrated into several downstream libraries, including \texttt{PyEPO}, \texttt{CVXPY}, and \texttt{cvxpylayers}, enabling MPAX to be used as a solver backend in end-to-end learning and decision-focused pipelines.

\end{itemize}

\section{Experiments}

In this section, we present computational results illustrating the four advanced features, namely across-hardware portability, batched solving, distributed optimization, and differentiable optimization, introduced in Section~\ref{sec:mpax}, along with benchmark comparisons against other implementations on standard benchmark sets.

\subsection{Across-hardware portability and hardware acceleration}

In this experiment, we evaluate the performance of \textsc{MPAX} across a diverse set of hardware platforms. By leveraging the XLA compiler, MPAX executes efficiently across this diverse hardware landscape. Importantly, the solver code and configuration remain identical across all runs; the only variable is the underlying hardware. This isolates hardware effects and ensures that differences in performance are not confounded by implementation or algorithmic changes.

\textbf{Problem instances.} We evaluate the performance of MPAX on two LP problems: the LP relaxations of a dense multi-dimensional knapsack instance and the sparse \texttt{dlr1} instance from the Mittelmann's benchmark set. These two cases are selected to assess \textsc{MPAX} under two distinct structural regimes, fully dense versus highly sparse, and to highlight how hardware characteristics interact with problem structure. Both instances are in the large-scale regime.

The LP relaxation of the multi-dimensional knapsack problem is formulated as follows.
\begin{equation}
\begin{aligned}
\min_{x\in[0,1]^m}\quad &  c^\top x \\
\text{s.t.}\quad &  Wx \le b .
\end{aligned}
\label{eq:knapsack}
\end{equation}
where $m$ is the number of items, $d$ is the number of knapsack dimensions, $c\in\mathbb{R}^m_+$ is the item-value vector, $W\in\mathbb{R}^{d\times m}_+$ is the weight matrix, and $b\in\mathbb{R}^d_+$ is the capacity vector.
The knapsack instance consists of 100,000 items and 1,000 dimensions with a fully dense constraint matrix. The \texttt{dlr1} instance is a sparse linear program featuring 9121907 variables, 1735470 constraints, and 18365107 nonzeros.

\textbf{Computational environment.} Our experiments utilize a wide range of modern hardware platforms, including an Intel Xeon Platinum 8462Y+ CPU (16 cores, 128 GB memory), NVIDIA GPUs (T4, L4, L40S, A100, H100, H200), and Google TPUs (v5e and v6e).  Detailed specifications for these platforms, including peak performance and memory bandwidth, are provided in Table \ref{tab:hardware-specs}. Note that because Google TPUs currently lack support for sparse matrix operations and double precision, they are evaluated only on the dense knapsack instance in single precision.

\begin{table}[h!]
\centering
\renewcommand{\arraystretch}{1.2}
\caption{Comparison of GPU/TPU compute capability, memory bandwidth, and memory capacity.}
\begin{tabular}{lccc}
\hline
\textbf{Device} & \textbf{Peak Performance (TFLOPS)$^\dagger$} & \textbf{Memory BW (GB/s)} & \textbf{Memory (GB)} \\ 
\hline
NVIDIA T4        & 8.1 / --          & 320   & 16  \\
NVIDIA L4        & 30.3 / --         & 300   & 24  \\
NVIDIA L40S      & 91.6 / --         & 864   & 48  \\
NVIDIA A100 PCIe & 19.5 / 9.7        & 1935  & 80  (HBM2e) \\
NVIDIA H100 SXM  & 66.9 / 33.5       & 3352  & 80  (HBM3) \\
NVIDIA H200 SXM  & 67 / 34           & 4800  & 141  (HBM3e) \\
Google TPU v5e   & 197 (BF16)        & 819   & 16  \\
Google TPU v6e   & 918 (BF16)        & 1600  & 32  \\
\hline
\end{tabular}
\label{tab:hardware-specs}
\footnotesize{\textit{$^\dagger$ NVIDIA numbers use FP32/FP64; TPU values are BF16 and not directly comparable.}}
\end{table}

\textbf{Results.}
Computational results are reported in Table~\ref{tab:hardware_runtime}. Since the algorithm is identical across different hardware platforms (up to cumulative machine-precision errors), we compare the runtime of the first 10{,}000 iterations of MPAX on different hardware by fixing the iteration limit to 10{,}000. We evaluate both single- and double-precision configurations. Both instances in all settings reached the iteration limit and the reported runtime is 10000 iterations.
Several observations in order: (i) On both CPUs and NVIDIA GPUs, single-precision runs are roughly twice as fast as their double-precision counterparts (with the exception of the CPU results for the dlr1 instance, discussed below). (ii) Relative to the high-end Intel 8462 CPU, the NVIDIA H200 GPU achieves a $32\times$ speedup in double precision and a $23\times$ speedup in single precision for the knapsack instance. The performance gap for dlr1 on CPUs can be attributed to the current JAX implementation: sparse operations on CPUs are not lowered to optimized libraries (e.g., MKL sparse routines) but are instead expressed via dense XLA operations, resulting in poor efficiency for large sparse problems. (iii) Even low-end GPUs such as the NVIDIA T4 and L4 outperform the high-end Intel CPU. (iv) Google TPUs deliver performance comparable to the NVIDIA L40S.

\begin{table}[t]
\centering
\caption{Runtime (in seconds) for the first 10{,}000 iterations of MPAX on two LP instances across different hardware platforms.}
\label{tab:hardware_runtime}
\begin{tabular}{lcccc}
\hline
\textbf{Hardware} 
& \multicolumn{2}{c}{\textbf{Knapsack}} 
& \multicolumn{2}{c}{\textbf{dlr1}} \\
\cline{2-3} \cline{4-5}
& Single precision & Double precision 
& Single precision & Double precision \\
\hline
Intel 8462Y+      & 58.68 & 133.06 & 1012.39 & 1243.18 \\
Nvidia T4         & 33.87 & 72.41  & 46.88  & 82.58  \\
Nvidia L4         & 32.33 & 65.48  & 43.75  & 79.75  \\
Nvidia L40S       & 11.47 & 22.96  & 14.40  & 27.94  \\
Nvidia A100       & 5.63  & 10.26  & 7.28   & 11.98  \\
Nvidia H100       & 3.05  & 5.50   & 4.45   & 7.00   \\
Nvidia H200       & 2.49  & 4.15   & 3.70   & 5.59   \\
Google TPU v5e   & 11.01 & --     & --     & --     \\
Google TPU v6e   & 7.47  & --     & --     & --     \\
\hline
\end{tabular}
\end{table}

\subsection{Batched solving}

We evaluate the batch-solving capability of \textsc{MPAX} using a two-stage stochastic cargo flight scheduling instance (the “storm” case) from~\cite{mulvey1995new}, where many second-stage LPs are supposed to be solved within each first-stage iteration.

\textbf{Problem instance.} This model originates from a classic air-cargo planning problem, where the first stage determines the number and types of aircraft assigned to a set of routes, subject to minimum flight frequencies and aircraft flow-balance constraints over the network. The second stage models cargo allocation under demand uncertainty, including shipped cargo, unmet demand, and unused capacity on each leg, with corresponding demand and flow-balance constraints.
Uncertainty enters only through the right-hand sides of 118 demand constraints. Following the construction of Mulvey and Ruszczyński (1995), we generate scenario-specific demands by scaling the nominal forecast by factors \{0.8, 0.9, 1.0, 1.1, 1.2\}, each with probability 0.2, independently across all demands. This results in a large-scale stochastic program with $5^{118}$ scenarios. The exact mathematical formulation of this problem is provided in Appendix~\ref{app:cargo}.

To evaluate batch-solving performance, we solve scenario subproblems using batch sizes of 1, 10, 100, 1{,}000, and 10{,}000, and compare the solve times of MPAX and Gurobi. Specifically, for each batch size, we measure the time required by each solver to solve the corresponding batch of instances in a single call.

\textbf{Computational environment.} We use an NVIDIA H100 SXM GPU with CUDA 12.4 for running MPAX and use Intel(R) Xeon(R) Platinum 8462Y+ @ 2.80 GHz with 128GB RAM and 16 threads for running Gurobi 13.0.0. For Gurobi, we solve the batched instances sequentially on the CPU. Each instance is solved with 16 threads enabled, so every Gurobi run fully utilizes the available CPU parallelism.

\begin{figure}[htb]
    \centering
    \includegraphics[width=0.7\linewidth]{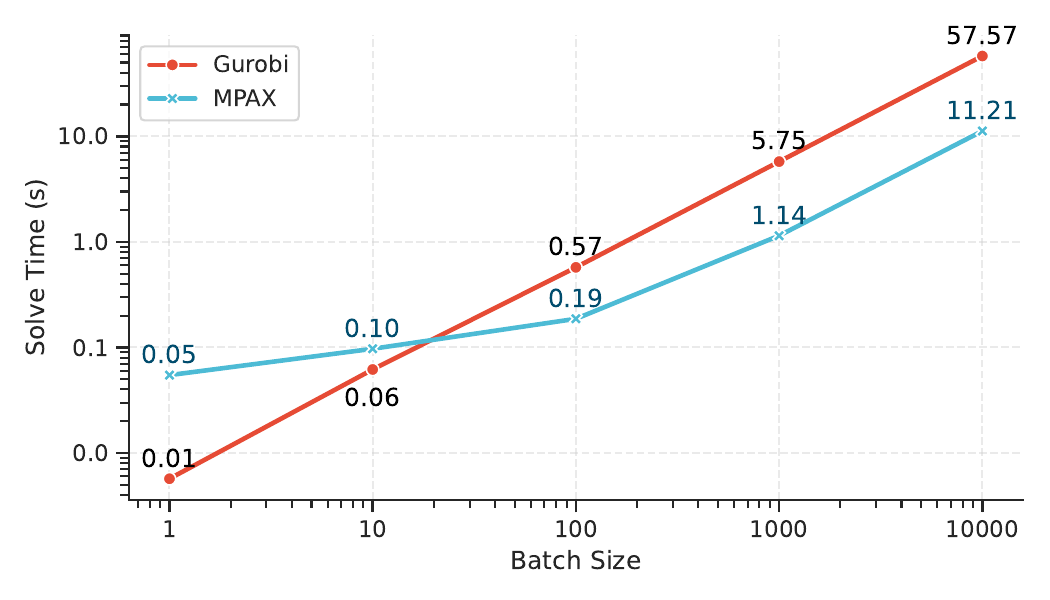}
    \caption{Solve time of MPAX and Gurobi for different batch sizes on the storm instance.}
    \label{fig:batch solve}
\end{figure}

\textbf{Results.}
The results are shown in Figure~\ref{fig:batch solve}. For a batch size of 1 (solving one LP at a time), Gurobi is about five times faster than MPAX (0.01 s vs. 0.05 s). However, as the batch size increases, MPAX exhibits much better scaling. At a batch size of 100, MPAX solves all problems in 0.19 s, while Gurobi requires 0.57 s. For the largest batch size of 10,000, the gap becomes substantial: MPAX takes 11.21 s, compared with Gurobi’s 57.57 s. Overall, for this small instance, the figure demonstrates that although Gurobi is faster for a single instance, MPAX scales dramatically better as the batch size grows.

\subsection{Distributed Optimization}

We evaluate \textsc{MPAX}’s distributed optimization capability on a large-scale multi-dimensional knapsack instance.

\textbf{Problem instance.}
We consider the knapsack model described in~\eqref{eq:knapsack} with 300{,}000 items and 3{,}000 dimensions, yielding an LP with approximately $9\times 10^8$ nonzeros. When multiple GPUs are used, the constraint matrix $A$ and its transpose $A^\top$ are partitioned by rows and distributed evenly across the devices to accelerate the dominant matrix–vector multiplications. The same instance is solved using 1--4 NVIDIA GPUs, and we report the wall-clock solve time.

\textbf{Computational environment.}
All experiments are conducted on a single node equipped with four NVIDIA H100 SXM GPUs with CUDA 12.4. The GPUs are interconnected via NVLink (NV6 links).

\begin{figure}[htb]
    \centering
    \includegraphics[width=0.5\linewidth]{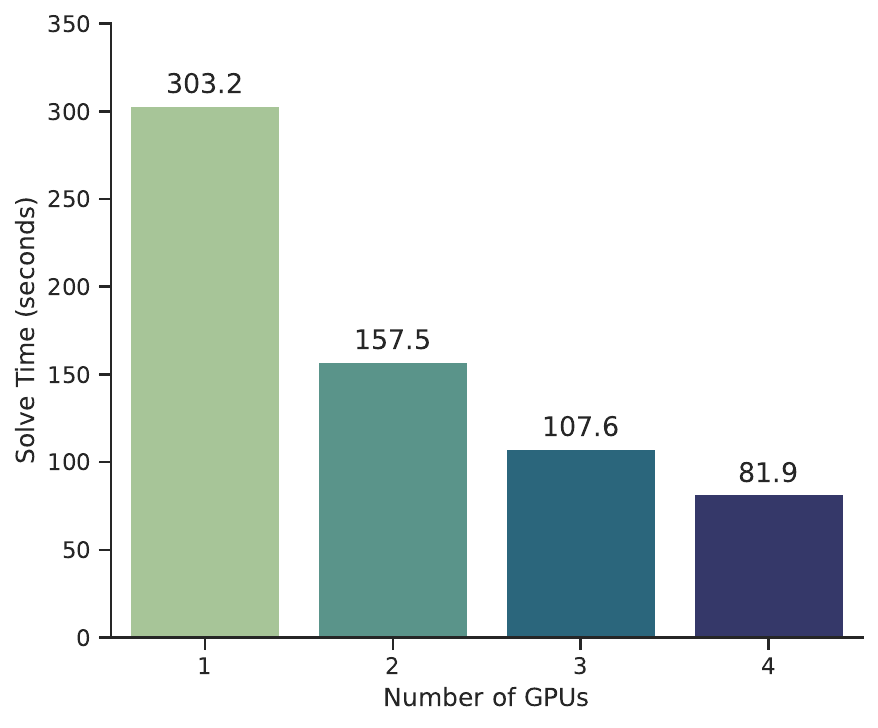}
    \caption{MPAX's distributed optimization performance on a multi-dimensional knapsack problem.}
    \label{fig:knapsack distributed optimization}
\end{figure}

\textbf{Results.} Across all configurations, MPAX terminates after 66{,}200 iterations; the measured solve time is shown in Figure~\ref{fig:knapsack distributed optimization}. As the number of GPUs increases, the solve time decreases from 303.2\,s (1 GPU) to 81.9\,s (4 GPUs), corresponding to a $3.7\times$ speedup. The scaling is close to linear, with mild diminishing returns at higher GPU counts, which is consistent with increased communication and synchronization overhead in distributed execution. These results demonstrate the effectiveness of distributed optimization with MPAX on a dense LP instance. We also note that a current limitation is that distributed optimization for sparse matrices requires more specialized sharding strategies to achieve favorable efficiency. Developing more sophisticated sharding schemes for sparse problems is left for future work.

\subsection{Differentiable Optimization}
\label{section:differentiable optimization}

Traditional solvers focus on finding an optimum of a mathematical program, whereas differentiable optimization focuses on computing the derivatives of the optimal solution to a mathematical program with respect to the parameters defining the program's objective and constraints. Differentiable optimization integrates optimization layers within machine learning models by enabling the computation of gradients through the solver~\cite{amos2017optnet,amos2018differentiable}. This capability is critical for end-to-end learning tasks where decisions made by the math programming problem depend on model predictions.

In this section, we present experiments evaluating the differentiable optimization capability of MPAX on an end-to-end Warcraft shortest path problem~\cite{poganvcic2020differentiation}. 

\begin{figure}[htb]
\begin{center}
\centerline{\includegraphics[width=0.75\columnwidth]{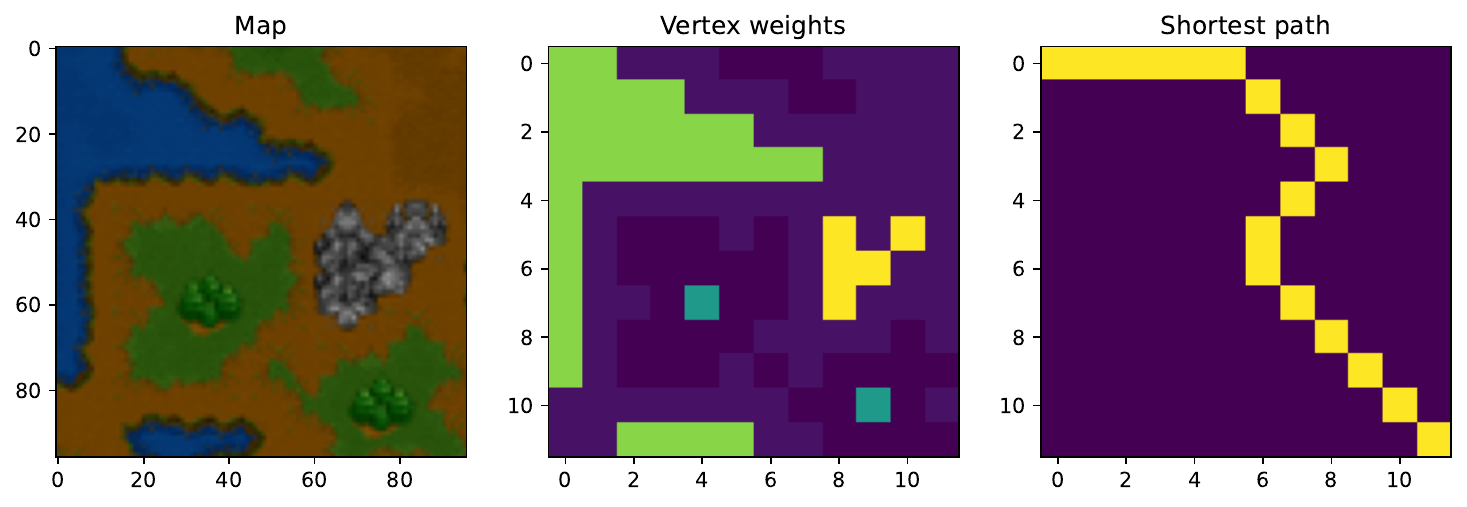}}
\caption{Warcraft shortest path dataset: Warcraft terrain map (left), vertex traversal cost (middle) and shortest path (right)}
\label{fig:warcraft example} 
\end{center}
\vskip -0.2in
\end{figure}

\textbf{Problem instance.}
We use the dataset introduced by \cite{poganvcic2020differentiation} for the Warcraft shortest path problem. An example of the Warcraft terrain map, vertex cost and shortest path in the dataset is provided in Figure \ref{fig:warcraft example}. In this task, the objective is to determine the shortest path between the top left and bottom right vertices on a Warcraft terrain map represented as a $k \times k$ 2D grid. The traversal cost of each vertex is unknown and dependent on the terrain type in the map image. The dataset includes Warcraft terrain maps of four different sizes, where $k \in \{12,18,24,30\}$. For each size, the training set contains 10,000 RGB images and the test set includes 1,000 RGB images.

In line with \cite{poganvcic2020differentiation}, we use the first five layers of ResNet \cite{he2016deep}, followed by a max-pooling operation, to extract the latent costs for the vertices. The LP formulation of the 2D grid shortest path problem is as follows. 
\begin{equation}
\label{eq:warcraft shortest path LP}
\begin{aligned}
\min & \quad \sum_{(i,j)\in E} c_{i,j} x_{i,j} \\
\text{s.t.} & \quad \sum_{j:(s,j) \in E} x_{s,j} - \sum_{i:(i,s) \in E} x_{i,s} = 1 \\
 & \quad \sum_{j:(t,j) \in E} x_{t,j} - \sum_{i:(i,t) \in E} x_{i,t} = -1 \\
 & \quad \sum_{j:(v,j) \in E} x_{v,j} - \sum_{i:(i,v) \in E} x_{i,v} = 0, \quad \forall v \in V \setminus \{s,t\} \\
& \quad 0 \leq x_{i,j} \leq 1, \quad \forall (i,j) \in E.
\end{aligned}
\end{equation}
where $V$ is the set of nodes, $s$ is the source node, $t$ is the sink node, $E$ is the set of edges, $c_{i,j}$ is the cost associated with travelling along edge $(i,j)$, and $x$ is the flow along edge $(i,j)$.
It is worth noting that the constraint matrix in this problem is very sparse. Therefore, we use the sparse matrix format when using MPAX.

\textbf{Computing environment.} 
We use an NVIDIA H100 SXM GPU, with CUDA 12.4, for running MPAX, FLAX, and PyTorch, and use Intel(R) Xeon(R) Platinum 8462Y+ @ 2.80 GHz with 128GB RAM and 16 threads for running Gurobi 13.0.0, OSQP 0.6.7, and SCS 3.2.7.

\textbf{Experiment details.} A differentiable optimization pipeline consists of a neural network component coupled with an optimization layer. We compare the following configurations: FLAX+MPAX, FLAX+OSQP, FLAX+SCS, and PyEPO+Gurobi. Here, FLAX is the neural network library in JAX, while PyEPO is a PyTorch-based differentiable optimization framework that integrates multiple optimization solvers. OSQP and SCS are accessed through the CvxpyQP wrapper in JAXOpt.

Since CvxpyQP does not allow modification of solver parameters for OSQP and SCS, we run both solvers with their default settings using 16 parallel CPU threads. Their default termination tolerances are $10^{-3}$ for OSQP and $10^{-4}$ for SCS.

For MPAX, we test termination tolerances $\epsilon \in \{10^{-3}, 10^{-4}, 10^{-5}, 10^{-6}\}$ while keeping all other options at their default values. When warm starting is enabled, the solution from the previous iteration is used as the initial point; otherwise, we initialize with an all-zero vector. Feasibility polishing is disabled during training. Since double precision is standard in linear programming, we evaluate MPAX exclusively using double-precision arithmetic to ensure consistency with common solver practices.

For Gurobi, we use the PyEPO package, which builds neural networks in PyTorch and solves linear programs using Gurobi. By default, when restricted to a single thread, Gurobi automatically selects among the primal simplex, dual simplex, and barrier methods. To ensure a fair comparison with the $\mathrm{r^2}$HPDHG algorithm used in MPAX, we explicitly enforce Gurobi’s barrier method without crossover by setting \texttt{Method=2} and \texttt{Crossover=0}. We also test different termination tolerances by setting \texttt{FeasibilityTol}, \texttt{OptimalityTol}, and \texttt{BarConvTol} (for the barrier method) to ${10^{-3}, 10^{-4}, 10^{-5}, 10^{-6}}$.

The model is trained over 10 epochs to minimize the SPO+ loss~\cite{elmachtoub2022smart} using the Adam optimizer with a learning rate of $10^{-4}$ and batch size of 70. We use the normalized regret as the evaluation metric. We set a time limit of 1 hour per epoch during training. All the experiment settings and metrics are the same as detailed in~\cite{tang2024pyepo}.

\begin{figure*}[htb]
\vskip 0.2in
\begin{center}
\centerline{\includegraphics[width=\textwidth]{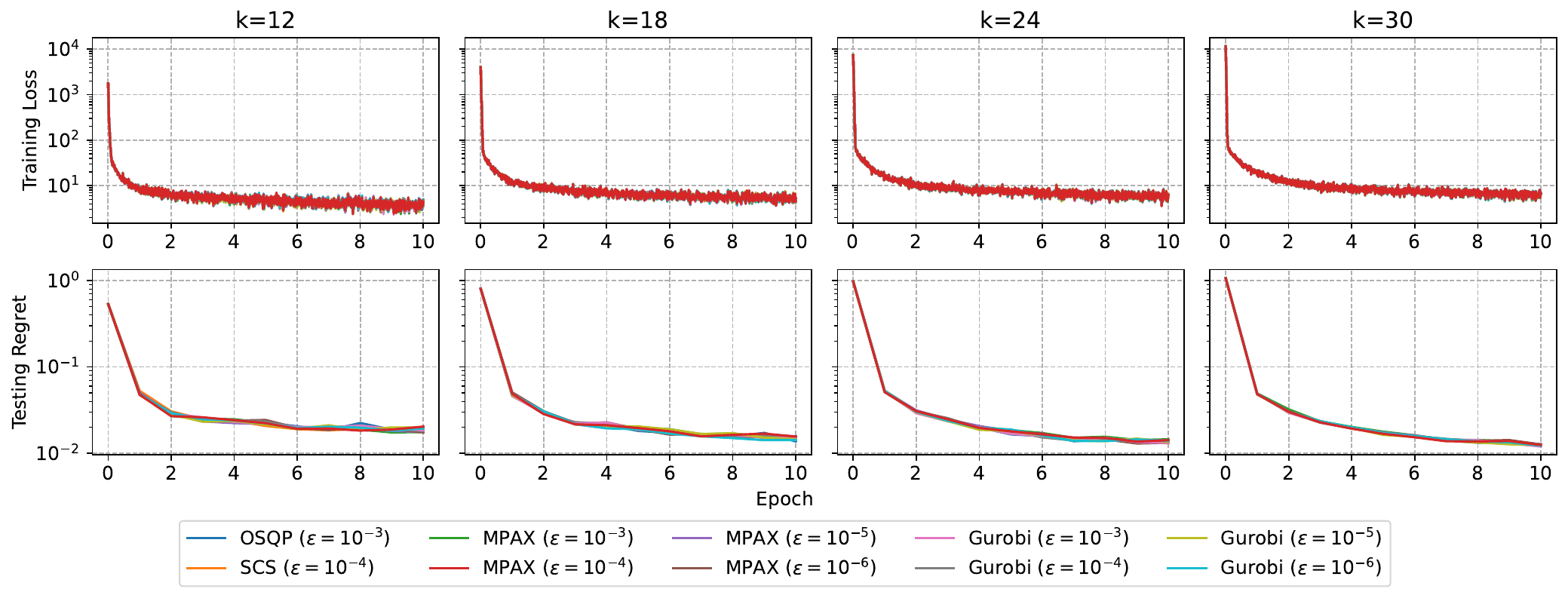}}
\caption{The loss and normalized regret curve of the Warcraft shortest path problem}
\label{fig: warcraft loss and regret}
\end{center}
\vskip -0.2in
\end{figure*}

\begin{table*}[htb]
\vskip 0.15in
\begin{center}
\begin{small}
\caption{Comparison of the training time per epoch and the average number of iterations for different strategies.}
\label{tab:warcraft time and iteration result}
\begin{tabular}{@{}cccrrrrrrrr@{}}
\toprule
\multirow{2}{*}{Framework} & \multirow{2}{*}{Algorithm} & \multirow{2}{*}{$\epsilon$} & \multicolumn{4}{c}{Training time per epoch} & \multicolumn{4}{c}{Average PDHG iteration} \\ \cmidrule(l){4-7} \cmidrule(l){8-11} 
 &  &  & k=12 & k=18 & k=24 & k=30 & k=12 & k=18 & k=24 & k=30 \\ \midrule
\multirow{8}{*}{FLAX + MPAX} & \multirow{4}{*}{$\mathrm{r^2HPDHG}$} & $10^{-3}$ & 9.9 & 20.0 & 23.0 & 38.1 & 736 & 1028 & 1625 & 2130 \\
 &  & $10^{-4}$ & 12.0 & 23.7 & 30.2 & 49.5 & 925 & 1302 & 1982 & 2536 \\
 &  & $10^{-5}$ & 13.0 & 25.0 & 31.1 & 52.9 & 1004 & 1413 & 2101 & 2683 \\
 &  & $10^{-6}$ & 13.8 & 26.2 & 34.4 & 56.5 & 1073 & 1512 & 2214 & 2808 \\
\cmidrule(l){2-11} 
 & \multirow{4}{*}{\begin{tabular}[c]{@{}c@{}}$\mathrm{r^2HPDHG}$\\ (warmstart)\end{tabular}} & $10^{-3}$ & 9.6 & 18.1 & 18.3 & 27.8 & 638 & 804 & 986 & 1172 \\
 &  & $10^{-4}$ & 11.3 & 21.5 & 21.5 & 35.1 & 813 & 1022 & 1252 & 1502 \\
 &  & $10^{-5}$ & 12.2 & 22.1 & 23.5 & 37.0 & 902 & 1130 & 1373 & 1634 \\
 &  & $10^{-6}$ & 13.2 & 23.9 & 24.9 & 40.7 & 983 & 1236 & 1485 & 1756 \\
\midrule
\multirow{1}{*}{FLAX + OSQP} & \multirow{1}{*}{ADMM} & $10^{-3}$ & 457.9 & 958.2 & 2320.7 & \multicolumn{1}{c}{-} & \multicolumn{1}{c}{-} & \multicolumn{1}{c}{-} & \multicolumn{1}{c}{-} & \multicolumn{1}{c}{-} \\
\midrule
\multirow{1}{*}{FLAX + SCS} & \multirow{1}{*}{ADMM} & $10^{-4}$ & 449.2 & 995.3 & 2291.7 & \multicolumn{1}{c}{-} & \multicolumn{1}{c}{-} & \multicolumn{1}{c}{-} & \multicolumn{1}{c}{-} & \multicolumn{1}{c}{-} \\
\midrule
\multirow{4}{*}{PyEPO + Gurobi} & \multirow{4}{*}{\begin{tabular}[c]{@{}c@{}}Barrier\\ (no crossover)\end{tabular}} & $10^{-3}$ & 14.1 & 37.5 & 88.0 & 170.4 & \multicolumn{1}{c}{-} & \multicolumn{1}{c}{-} & \multicolumn{1}{c}{-} & \multicolumn{1}{c}{-} \\
 &  & $10^{-4}$ & 16.0 & 41.9 & 97.3 & 194.8 & \multicolumn{1}{c}{-} & \multicolumn{1}{c}{-} & \multicolumn{1}{c}{-} & \multicolumn{1}{c}{-} \\
 &  & $10^{-5}$ & 13.1 & 34.0 & 80.7 & 170.6 & \multicolumn{1}{c}{-} & \multicolumn{1}{c}{-} & \multicolumn{1}{c}{-} & \multicolumn{1}{c}{-} \\
 &  & $10^{-6}$ & 16.7 & 42.9 & 95.2 & 179.2 & \multicolumn{1}{c}{-} & \multicolumn{1}{c}{-} & \multicolumn{1}{c}{-} & \multicolumn{1}{c}{-} \\
\bottomrule
\end{tabular}
\end{small}
\end{center}
\vskip -0.1in
\end{table*}

\textbf{Results.}
Table~\ref{tab:warcraft time and iteration result} reports the training time per epoch and the average number of PDHG iterations required for MPAX to converge. A few observations in order: (i) both FLAX+OSQP and FLAX+SCS become prohibitively slow as the grid size increases: when $k=30$, neither solver is able to complete a single training epoch within the one-hour time limit. (ii) Across all configurations, FLAX+MPAX consistently achieves the strongest performance. Its advantage over FLAX+JAXOpt (OSQP/SCS) and PyEPO+Gurobi becomes increasingly pronounced as the problem size grows. For example, at $k=30$, FLAX+MPAX completes training in roughly one third of the time required by PyEPO+Gurobi. (iii) Within MPAX, warm starting further improves efficiency: it reduces the number of PDHG iterations by approximately 35\%, leading to a corresponding reduction in training time of about 25\%. These results highlight both the scalability advantages of MPAX and the effectiveness of warm starting for large-scale differentiable optimization tasks.

Figure~\ref{fig: warcraft loss and regret} presents the training loss and testing regret as functions of training epochs for different solver configurations and termination tolerances. Two observations are noteworthy. First, for this differentiable optimization task, a tolerance of $10^{-3}$ is sufficient to achieve comparable performance in terms of both training loss decay and testing regret reduction. This indicates that the model’s performance is relatively insensitive to the accuracy of the underlying optimization solver. Second, since the learning curves across solvers are essentially identical when measured in epochs, computational efficiency is primarily determined by the time required to complete each epoch, as quantified in Table~\ref{tab:warcraft time and iteration result}.

\subsection{Benchmark Results}

In this section, we compare the numerical performance of MPAX with cuPDLPx on LP benchmarks and with PDQP.jl on QP benchmarks, focusing on single-instance solves on standard benchmark sets. The goal of the experiment is not to do a comprehensive study compared with all other solvers, but to compare the effectiveness of MPAX with other implementations of the same algorithms to get more insights. Note the advance features that are mentioned above are unique at MPAX.

In this section, we compare the numerical performance of MPAX with cuPDLPx on LP benchmarks and with PDQP.jl on QP benchmarks, focusing on single-instance solves over standard benchmark sets. The goal of these experiments is not to provide a comprehensive comparison against all available solvers, but rather to assess the effectiveness of MPAX relative to other implementations of the {same underlying algorithms}, thereby isolating the impact of the software and system design choices. We note that several advanced features discussed earlier are unique to MPAX and are therefore not reflected in these comparisons.

\textbf{Problem instances.}
For LP, we use  two widely used LP benchmark datasets: MIP Relaxations, which consists of 383 instances derived from root-node LP relaxations
of MIPLIB 2017, and Mittelmann’s LP benchmark dataset, which contains 49 LP instances.
For QP, we use 134 convex QP instances from the Maros–Meszaros test set and 33 convex QP instances from
QPLIB. 

\textbf{Computing environment.} We use an NVIDIA H100 SXM GPU, with CUDA 12.4. The experiments of MPAX are performed in JAX 0.7.2, and the experiments of PDQP are performed in Julia 1.11.5.

\textbf{Implementation details.} MPAX and cuPDLPx implement the same algorithm for LPs, while MPAX and PDQP.jl implement the same algorithm for QPs. The primary differences only lie in the programming environments. cuPDLPx is written in C, MPAX in JAX, and PDQP.jl in Julia. Because both JAX and Julia rely on just-in-time compilation, we exclude compilation time and report only execution time to ensure a fair comparison.

MPAX uses the relative KKT error for termination, which checks primal feasibility, dual feasibility, and primal-dual gap. This is the same as both cuPDLPx and PDQP.jl.
For LP instances, we set the time limit to 15000 seconds and set $\epsilon = 10^{-4}$ as the termination criterion.
For QP instances, we set the time limit to 3600 seconds and set $\epsilon = 10^{-3}$ for moderate accuracy and $ \epsilon = 10^{-6}$ for high accuracy.
We report the shifted geometric mean of the solve time required across instances for each solver. More precisely, the shifted geometric mean is defined as $\left( \prod_{i=1}^n (t_i + \Delta) \right)^{1/n} - \Delta$ and we shift by $\Delta = 10$.  If the instance cannot be solved within the time limit, we set the solve time as the time limit. We also report Count, the number of instances that terminate with a conclusive status (optimal, infeasible, or unbounded) within the time limit.

\textbf{Results.}
The LP and QP benchmark results are summarized in Tables~\ref{tab:lp benchmark} and~\ref{tab:qp benchmark}. For LPs, cuPDLPx achieves a 1.05$\times$ speedup over \textsc{MPAX} on the MIPLIB set and a 1.18$\times$ speedup on the Mittelmann set. This difference is expected, as cuPDLPx is implemented in high-performance C and benefits from a highly optimized sparse linear algebra backend.

For QPs, \textsc{MPAX} outperforms PDQP.jl by approximately 1.25$\times$–1.6$\times$ across the benchmark instances. We attribute this gap primarily to implementation maturity: the current PDQP.jl codebase was developed in 2023 and has not yet been updated or tuned for the most recent CUDA toolchain and GPU architectures.

Overall, despite being written in Python, \textsc{MPAX} delivers performance competitive with highly optimized C implementations. This highlights the effectiveness of XLA compilation in transforming high-level JAX code into efficient accelerator-executable kernels.

\begin{table}[htbp]
\centering
\caption{Solve time (SGM10) in seconds of MPAX and cuPDLPx on LP benchmarks.}
\begin{tabular}{lcccc}
\toprule
Dataset & Tolerance & Solver & Count & SGM10 Time \\
\midrule
\multirow{2}{*}{Mittelmann}     & \multirow{2}{*}{$10^{-4}$} & MPAX & 48 & 22.61 \\
                                &        & cuPDLPx & 48 & 19.11 \\
\midrule
\multirow{2}{*}{MIPLIB}         & \multirow{2}{*}{$10^{-4}$} & MPAX & 382 & 4.63 \\
                                &          & cuPDLPx & 382 & 4.40 \\
\bottomrule
\end{tabular}
\label{tab:lp benchmark}
\end{table}

\begin{table}[htbp]
\centering
\caption{Solve time (SGM10) in seconds of MPAX and PDQP.jl on QP benchmarks.}
\begin{tabular}{lcccc}
\toprule
Dataset & Tolerance & Solver & Count & SGM10 Time \\
\midrule
\multirow{2}{*}{Maros-Meszaros} & \multirow{2}{*}{$10^{-3}$} & MPAX & 130 & 5.29 \\
                                &          & PDQP.jl & 130 & 8.50 \\
\midrule
\multirow{2}{*}{Maros-Meszaros} & \multirow{2}{*}{$10^{-6}$} & MPAX & 122 & 23.73 \\
                                &          & PDQP.jl & 121 & 35.41 \\
\midrule
\multirow{2}{*}{QPLIB}          & \multirow{2}{*}{$10^{-3}$} & MPAX & 31  & 19.44 \\
                                &          & PDQP.jl & 29  & 24.28 \\
\midrule
\multirow{2}{*}{QPLIB}          & \multirow{2}{*}{$10^{-6}$} & MPAX & 27  & 36.92 \\
                                &          & PDQP.jl & 28  & 52.88 \\
\bottomrule
\end{tabular}
\label{tab:qp benchmark}
\end{table}

\section{Conclusion and Future Work}

In this paper, we present MPAX, a hardware-accelerated, batchable, distributable, and differentiable solver for linear programming and quadratic programming by exploiting modern machine learning infrastructure JAX. Benefiting from the powerful JAX transformations, MPAX provides advanced features, including across-hardware portability, batched solving, distributed optimization, and automatic differentiation. Our numerical results show that MPAX achieves substantial acceleration on accelerators, competitive benchmark performance against GPU-based solver codebases, and strong scalability in batched and multi-GPU settings. While MPAX currently supports linear and quadratic programming, future work will extend its capability to more general mathematical programming problems and specialized modules for common machine-learning tasks.

\section*{Acknowledgement}
Haihao Lu is supported by AFOSR Grant No. FA9550-24-1-0051 and ONR Grant No. N000142412735. Zedong Peng is supported by ONR Grant No. N000142412735. Jinwen Yang is supported by AFOSR Grant No. FA9550-24-1-0051.

\bibliographystyle{amsplain}
\bibliography{ref-mpax}

\newpage
\appendix
\section{Solver Options}
\label{appendix:solver options}
\begin{itemize}
    \item \textbf{General Options}
\begin{longtable}{|c|c|p{9.5cm}|}
\hline
\textbf{Option} & \textbf{Default} & \textbf{Description} \\ \hline
\texttt{verbose} & \texttt{False} & Whether to print solver output to stdout. \\ \hline
\texttt{debug} & \texttt{False} & Whether to print additional debugging information to stdout. \\ \hline
\texttt{jit} & \texttt{True} & Whether to enables Just-In-Time (JIT) compilation. \\ \hline
\texttt{unroll} & \texttt{False} & Whether to unroll iteration loops. \\ \hline
\texttt{display\_frequency} & \texttt{10} & Frequency (in every termination check) to print solver output. \\ \hline
\texttt{warm\_start} & \texttt{False} & Whether to perform warm starting. \\ \hline
\texttt{feasibility\_polishing} & \texttt{False} & Whether to perform feasibility polishing. \\ \hline
\end{longtable}

\item \textbf{Termination Parameters}
\begin{longtable}{|c|c|p{9cm}|}
\hline
\textbf{Option} & \textbf{Default} & \textbf{Description} \\ \hline
\texttt{eps\_abs} & \texttt{1e-4} & Absolute tolerance for convergence\\ \hline
\texttt{eps\_rel} & \texttt{1e-4} & Relative tolerance for convergence\\ \hline
\texttt{eps\_primal\_infeasible} & \texttt{1e-8} & Primal infeasibility tolerance\\ \hline
\texttt{eps\_dual\_infeasible} & \texttt{1e-8} & Dual infeasibility tolerance\\ \hline
\texttt{eps\_feas\_polish} & \texttt{1e-6} & Tolerance for feasibility polishing\\ \hline
\texttt{iteration\_limit} & \texttt{max\_int} & Maximum number of iterations to run\\ \hline
\end{longtable}

\end{itemize}

\section{Stochastic cargo flight scheduling problem}\label{app:cargo}

We consider a two-stage stochastic cargo flight scheduling problem. The first-stage decisions specify the number of sorties assigned to each candidate route for each aircraft type. After the demand scenario is realized, second-stage decisions determine how cargo is delivered along the scheduled routes, possibly with transshipment, and quantify any unmet demand. The goal is to minimize the expected total cost, including operating costs, cargo handling costs, and penalties for undelivered cargo, subject to payload capacity, landing limits, service requirements, and flying-hour constraints.

\paragraph{Sets and Parameters}
\begin{align*}
i &= \text{scenario index, } i=1,\ldots,S;\\
p_i &= \text{probability of scenario } i;\\
\mathcal{N} &= \text{the set of nodes;}\\
b_i(m,n) &= \text{the amount of cargo to be shipped from nodes $m$ to $n$ in scenario $i$;}\\
f(m,n) &= \text{the minimum number of flights from $m$ to $n$;}\\
\Pi &= \text{the set of routes, i.e., selected sequences of nodes $n_1,n_2,\ldots,n_l$;}\\
l(\pi) &= \text{the number of nodes in route $\pi$;}\\
v(\pi,j) &= \text{the $j$th node in route $\pi$, $j=1,\ldots,l(\pi)$;}\\
U(m,n) &= \{\pi\in\Pi: m=v(\pi,j_1),\ n=v(\pi,j_2),\ j_1<j_2\};\\
V_1(n) &= \{\pi\in\Pi: n=v(\pi,1)\};\\
V_l(n) &= \{\pi\in\Pi: n=v(\pi,l(\pi))\};\\
W(n) &= \{\pi\in\Pi: n=v(\pi,j)\text{ for some } j\};\\
\sigma(n) &= \text{the maximum number of landings allowed in node $n$;}\\
\mathcal{A} &= \text{the set of aircraft types;}\\
d(a) &= \text{the maximum payload of an aircraft of type $a$;}\\
h^{\max}(a) &= \text{the maximum flying hours of aircrafts of type $a$;}\\
h^{\min}(a) &= \text{the minimum flying hours of aircrafts of type $a$;}\\
h(\pi,a) &= \text{the flying hours necessary for an aircraft of type $a$ to pass route $\pi$;}\\
c(a) &= \text{the cost of an operating hour of an aircraft of type $a$;}\\
q &= \text{the unit cost of cargo handling;}\\
\rho &= \text{the penalty for each undelivered cargo unit.}
\end{align*}

\paragraph{Decision variables.}

\begin{align*}
x(\pi,a) &=
\begin{aligned}[t]
&\text{the number of sorties on route $\pi$ by aircraft type $a$;}
\end{aligned}\\
d_i(\pi,m,n) &=
\begin{aligned}[t]
&\text{the amount of cargo delivered directly from nodes $m$ to $n$ over route $\pi$}\\
&\text{in scenario $i$;}
\end{aligned}\\
t_i(\pi,m,k,n) &=
\begin{aligned}[t]
&\text{the amount of cargo from nodes $m$ to $n$ to be moved to the transshipment}\\
&\text{node $k$ over route $\pi$ in scenario $i$;}
\end{aligned}\\
s_i(\pi,k,n) &=
\begin{aligned}[t]
&\text{the amount of transshipment cargo to be moved from the transshipment}\\
&\text{node $k$ to node $n$ over route $\pi$ in scenario $i$;}
\end{aligned}\\
y_i(m,n) &=
\begin{aligned}[t]
&\text{the amount of undelivered cargo from nodes $m$ to $n$ in scenario $i$;}
\end{aligned}\\
z_i(\pi,j) &=
\begin{aligned}[t]
&\text{the unused capacity on leg $j$ of route $\pi$, $j=1,\ldots,l(\pi)-1$ in scenario $i$.}
\end{aligned}
\end{align*}

\paragraph{Objective function}

\begin{align}
\min\ &\Bigg[\sum_{\pi\in\Pi}\sum_{a\in\mathcal{A}} c(a)h(\pi,a)x(\pi,a) \nonumber\\
&\quad + \sum_{i=1}^{S} p_i\Bigg(q\sum_{\pi\in\Pi}\sum_{(m,n)\in I(\pi)}\Big(d_i(\pi,m,n)+s_i(\pi,m,n)+\sum_{k\in\mathcal{N}} t_i(\pi,m,n,k)\Big)  + \rho\sum_{m\in\mathcal{N}}\sum_{n\in\mathcal{N}} y_i(m,n)\Bigg)\Bigg]. \label{eq:obj}
\end{align}

\paragraph{Cargo balance constraint}
\begin{equation}
    \sum_{\pi\in\Pi}\Big(d_i(\pi,m,n)+\sum_{k\in\mathcal{N}} t_i(\pi,m,k,n)\Big)+y_i(m,n)
\ge b_i(m,n), \qquad \forall m,n\in\mathcal{N},\ i=1,\ldots,S \label{eq:cargo1}
\end{equation}

\begin{equation}
    \sum_{\pi\in\Pi}\sum_{m\in\mathcal{N}} t_i(\pi,m,k,n)
= \sum_{\pi\in\Pi} s_i(\pi,k,n), \qquad \forall k,n\in\mathcal{N},\ i=1,\ldots,S
\label{eq:cargo2}
\end{equation}

\paragraph{Payload balance constraints}

\begin{equation}
\begin{aligned}
\sum_{k\in\mathcal{N}}\Big(d_i(\pi,v(\pi,1),k)+s_i(\pi,v(\pi,1),k)+\sum_{n\in\mathcal{N}} t_i(\pi,v(\pi,1),k,n)\Big) = \sum_{a\in\mathcal{A}} d(a)x(\pi,a)-z_i(\pi,1),\\
\ \forall \pi\in\Pi,\ i=1,\ldots,S 
\label{eq:payload1}
\end{aligned}
\end{equation}

\begin{equation}
\begin{aligned}
\sum_{k\in\mathcal{N}}\Big(d_i(\pi,v(\pi,j),k)+s_i(\pi,v(\pi,j),k)+\sum_{n\in\mathcal{N}} t_i(\pi,v(\pi,j),k,n)\Big)\\
-\sum_{k\in\mathcal{N}}\Big(d_i(\pi,k,v(\pi,j))+s_i(\pi,k,v(\pi,j))+\sum_{n\in\mathcal{N}} t_i(\pi,k,v(\pi,j),n)\Big)\\
= z_i(\pi,j-1)-z_i(\pi,j), \qquad j=2,\ldots,l(\pi)-1,\ \forall \pi\in\Pi,\ i=1,\ldots,S
\label{eq:payload2}
\end{aligned}
\end{equation}

\paragraph{Frequency constraint}
\begin{align}
\sum_{a\in\mathcal{A}}\sum_{\pi\in U(m,n)} x(\pi,a) \ge f(m,n), \qquad \forall m,n\in\mathcal{N}, \label{eq:freq}
\end{align}

\paragraph{Landing and take-off balance constraint}
\begin{align}
\sum_{\pi\in V_1(n)} x(\pi,a) = \sum_{\pi\in V_l(n)} x(\pi,a), \qquad \forall a\in\mathcal{A},\ n\in\mathcal{N}, \label{eq:landtake}\\
\sum_{a\in\mathcal{A}}\sum_{\pi\in W(n)} x(\pi,a) \le \sigma(n), \qquad \forall n\in\mathcal{N}, \label{eq:landlimit}
\end{align}

\paragraph{Flying hours balance constraint}
\begin{align}
h^{\min}(a) \le \sum_{\pi\in\Pi} x(\pi,a)h(\pi,a) \le h^{\max}(a), \qquad \forall a\in\mathcal{A}. \label{eq:hours}
\end{align}

\end{document}